\documentclass{amsart}

\usepackage{fullpage}

\usepackage{amsmath}
\usepackage{amsfonts}
\usepackage{amsthm}

\newtheorem{prop}{Proposition}

\newcommand\SP{S$\cdot$P}

\title{Improved Bounds on the Sizes of \SP\ Numbers}
\author{Paul Myer Kominers and Scott Duke Kominers}
\thanks{The second author gratefully acknowledges the support of a Harvard Mathematics Department Highbridge Fellowship.}
\address{\newline\indent Student, Department of Mathematics, Massachusetts Institute of Technology}
\email{pkoms@mit.edu}
\address{\newline\indent Student, Department of Mathematics, Harvard University\newline\indent c/o 8520
  Burning Tree Road\newline \indent Bethesda, MD 20817}
\email{kominers@fas.harvard.edu}

\begin{document}
\maketitle
\section{Introduction}
A number which is \emph{\SP\ in base $r$} is a positive integer
which is equal to the sum of its base-$r$ digits multiplied by the product of its base-$r$ digits.
That is, $a_nr^n+\cdots + a_1r+a_0$ (here and hereafter, $0\leq a_i<r$ for all $0\leq i\leq n$) is \SP\ if and only if $$a_nr^n+\cdots + a_1r+a_0=a_0\cdots a_n(a_0+\cdots +a_n).$$  For example, $144=1\cdot 4\cdot 4\cdot(1+4+4)$ is \SP\ in base $10$ and $6=1\cdot 2\cdot (1+2)$ is \SP\ in base $4$.

Parame\'swaran \cite{P} conjectured that the number of base-$10$ \SP\ numbers is finite.  Several authors subsequently gave proofs of Parame\'swaran's conjecture and generalizations to other bases (see \cite{L}), as well as
enumerations of \SP\ numbers (see \cite{e1,e2}).  Recently, Shah Ali \cite{SA} gave a new argument proving that the number of base-$r$ \SP\ numbers is finite for any $r>1$.
In his proof, Shah Ali \cite{SA} obtained the first effective bound on the sizes of \SP\ numbers:
\begin{prop}(\cite{SA})\label{old}
A number which is \SP\ in base $r>1$ has at most $2r(r-1)^2$ digits.
\end{prop}
However, a quick check in the case $r=2$ shows that this bound is far
from sharp.  While Proposition \ref{old} shows that a base-$2$ \SP\
number can have at most $4$ digits, quick analysis shows that there is a
unique base-$2$ \SP\ number, $1$.  Indeed, if $a_n2^n+\cdots +
a_12+a_0$ is \SP\ in base $2$ then $a_i=1$ for $0\leq i\leq n$.
However, we then must have $$2^{n+1}-1=2^n+\cdots + 2^0=a_n2^n+\cdots
+a_0=a_n+\cdots +a_0=1+\cdots + 1=n+1;$$ it follows easily that $n=0$.

\section{A Sharp Bound}
Modifying Shah Ali's \cite{SA} method, we obtain an improved bound on the number of digits in a base-$r$ \SP\ number. As we will discuss in Section \ref{remarks}, our bound is sharp in the case $r=2$. 
\begin{prop}\label{new}
A number which is \SP\ in base $r>1$ has at most $2(r-1)^3-2(r-1)+1=2(r-1)(r^2-2r)+1$ digits.
\end{prop}
\begin{proof}
Let $a_nr^n+\cdots + a_1r+a_0$ be \SP\ in base $r$ with $n\geq 0$, so that
$$a_nr^n+\cdots + a_1r+a_0=a_0\cdots a_n(a_0+\cdots +a_n).$$  Then, $0<
a_i<r$ for all $0\leq i\leq n$, so that we have
\begin{equation}\label{new!}\left(\min_{0\leq i \leq n}\{a_i\}\right)\cdot \left(\frac{r^{n+1}-1}{r-1}\right)\leq a_0\cdots a_n(a_0+\cdots +a_n).\end{equation}
Then, since $\min_{0\leq i \leq n}\{a_i\}>0$, we may divide both sides of (\ref{new!}) by $\min_{0\leq i \leq n}\{a_i\}$ and obtain
\begin{align}
\label{hip}\frac{r^{n+1}-1}{r-1}\leq \frac{a_0\cdots a_n}{\min_{0\leq i \leq n}\{a_i\}}(a_0+\cdots +a_n)\leq (r-1)^n(n+1)(r-1),\end{align} as $\max_{0\leq i \leq n}\{a_i\}<r$.  Rearranging (\ref{hip}) gives
\begin{equation}\label{eq1}\left(1+\frac{1}{r-1}\right)^{n+1}=\left(\frac{r}{r-1}\right)^{n+1}\leq (n+1)(r-1)+\frac{1}{(r-1)^{n+1}}.\end{equation}
By the Binomial Theorem, the left side of (\ref{eq1}) is equal to
$$1+\frac{n+1}{r-1}+\frac{(n+1)n}{2(r-1)^2}+\cdots,$$ hence we obtain
\begin{equation}\label{eq2}1+\frac{n+1}{r-1}+\frac{(n+1)n}{2(r-1)^2}\leq (n+1)(r-1)+\frac{1}{(r-1)^{n+1}}.\end{equation}
Simplifying (\ref{eq2}), we find
\begin{equation}\label{eq3}n\leq 2(r-1)^3-2(r-1)+\frac{2}{n+1}\left(\frac{1}{(r-1)^{n-1}}-(r-1)^2\right).\end{equation}
Now, since $r>1$ and $n\geq 0$, we have that $(r-1)^{1-n}-(r-1)^2\leq 0$.  It then follows from (\ref{eq3}) that $n\leq 2(r-1)^3-2(r-1)$.
\end{proof}

\section{Remarks}\label{remarks}
Use of a computer algebra system suggests that \cite{L} implies the effective bound \begin{equation}\label{conk}n+1\leq \frac{W_{-1}\left(\frac{\log (r-1)-\log (r)}{r}\right)}{\log (r-1)-\log(r)}\end{equation} on the number $(n+1)$ of digits of a base-$r$ \SP\ number.  Here, $W_{-1}(\cdot)$ is the $(-1)^{\text{st}}$ analytic branch of the \emph{Lambert $W$-Function}, the multivalued inverse of the function $f(W)=We^W$.  (Weisstein \cite{Wei} summarizes the fundamental properties of the $W$-function. Corless \textit{et al.} \cite{Cea} survey several relevant applications and present an efficient method of evaluating the $W$-function to arbitrary precision.) Although we have been unable to prove the bound (\ref{conk}), we have verified it for $1<r\leq 999$.  

The bound given in Proposition \ref{new} is sharp for the case $r=2$; this is a 75\% improvement on the bound of Shah Ali's \cite{SA} Proposition \ref{old}.  Furthermore, although (\ref{conk}) is generally far smaller than the cubic bound of Proposition \ref{new}, (\ref{conk}) gives at best that an \SP\ number in base $2$ has no more than two digits.  Thus, our Proposition \ref{new} is the first sharp bound found for the case $r=2$.

For the case $r=10$, Proposition \ref{new} gives a bound of $1441$ digits, an 11\% improvement upon Proposition \ref{old}.  However, this bound is far weaker than the bounds given in \cite{L}, which show that a base-$10$ \SP\ number can have at most $60$ digits.

\section*{Acknowledgements}
The authors greatly appreciate the referee's helpful comments and suggestions on an earlier draft of the work.


\begin{thebibliography}{9}

\bibitem{L}Paul Belcher, H. J. Godwin, Andrew Lobb, Nick Lord, K. Robin McLean and Phillip Williams, On \SP\ numbers, \textit{Math. Gaz.} \textbf{82} (March 1998) 
pp.\ 72--75.
\bibitem{e2}Ezra Bussmann, \SP\ numbers in bases other than 10, \textit{Math. Gaz.} \textbf{85} (July 2001) pp. 245--248.
\bibitem{Cea}R.\ M.\ Corless, G.\ H.\ Gonnet, D.\ E.\ G.\ Hare, D.\ J.\ Jeffrey, and D.\ E.\ Knuth, On the Lambert $W$ function, \textit{Adv.\ Comp.\ Math.}\ \textbf{5} (December 1996) pp.\ 329-359.
\bibitem{e1}K. Robin McLean, There are only three \SP\ numbers,
  \textit{Math. Gaz.} \textbf{83}
(March 1999) pp. 32--39.
\bibitem{P}S. Parame\'swaran, Numbers and their digits -- a structural pattern, \textit{Math. Gaz.} \textbf{81} (July 1997) p. 263.
\bibitem{SA}H. A. Shah Ali, The number of \SP\ numbers is finite, \textit{Math. Gaz.} \textbf{92} (March 2008) pp. 64--65.
\bibitem{Wei}Eric W. Weisstein, Lambert's $W$-function, \textit{CRC Concise Encyclopedia of Mathematics}, CRC Press (2003) pp. 1684--1685.
\end{thebibliography}
\end{document}